\documentclass[12pt]{article}
\usepackage{amsmath,amsthm,amssymb,amscd,}
 \theoremstyle{definition}
 
 \theoremstyle{remark}

 \numberwithin{equation}{section}

\title
{Addendum to Maximal   regularity  and Hardy spaces
}
\author{ Pascal Auscher, Fr\'ed\'eric Bernicot and Jiman Zhao
  \\
 \small Universit\'e de Paris-Sud, Orsay et CNRS 8628, 91405 Orsay Cedex, France
\\
\small {\em E-mail address:} {pascal.auscher@math.u-psud.fr}\\
\small Universit\'e de Paris-Sud, Orsay et CNRS 8628, 91405 Orsay Cedex, France
\\
\small {\em E-mail address:} {Frederic.Bernicot@math.u-psud.fr}\\
\small School of Mathematical Sciences,
Beijing Normal University \\
\small Beijing 100875,
P.R. China\\
\small {\em E-mail address:} {jzhao@bnu.edu.cn}\\ 
}
\date{October 16, 2008}
\begin{document}
\maketitle

\begin{abstract} We correct an inaccuracy in \cite{ABZ}.

\vspace{0.1in}

\vspace{0.1in}
\end{abstract}

Proposition 3.14 in \cite{ABZ}  states as follows: Let $H_1$, $H_2$ be any  Hardy spaces obtained via an atomic decomposition. Any linear operator mapping the atoms of $H_1$ into a bounded set in $H_2$ has a bounded extension
from $H_1$ into $H_2$.

 M. Bownik \cite{B}, based on an example of Y.Meyer in Meyer, Taibleson and Weiss \cite{MTW}, showed that the  Hardy space norm and the finite  Hardy space norm may not be equivalent on the finite Hardy space (defined by restricting atomic decompositions to be finite sums and the atoms are $L^\infty$-atoms in the sense of Coifman-Weiss).  Hence, Proposition 3.14 is not correct as stated.

The 
 paper by Meda, Sj\"ogren and Vallarino \cite{MSV} establishes, among other things, that if one replaces $L^\infty$ atoms by $L^2$ atoms, the equivalence holds. Hence, for an operator to have a bounded extension  it suffices it is bounded on $L^2$-atoms. Moreover, the extension coincides with the original operator on $H^1\cap L^2$.  So Proposition 3.14  is correct if $H_{1}$ is the original Coifman-Weiss Hardy space and atoms are $L^2$-atoms.  As the atoms in \cite{ABZ} are $L^2$-atoms on a space of homogeneous type, this applies directly to the spaces $H^1_{z}(X)$ and (with little extra work)   $H^1_{r}(X)$ defined  in \cite{ABZ}: the maximal regularity operator  and its adjoint have  the boundedness property  announced in Theorem 2.1 there. 
 
 Another possibility to reach our conclusion is to apply the weak-type (1,1) result of Coulhon and Duong \cite{CD} under our hypotheses. This implies that the maximal regularity operator and its adjoint  are defined on $L^1$, thus on the Hardy space. Taking approximations of the these operators by truncating integrals in their definitions, one can adapt the proof of Lemma 1 in Chapter VI of Meyer's book \cite{M}.


\begin{thebibliography}{9999}


\bibitem{ABZ} Auscher, Pascal; Bernicot, Fr\'ed\'eric; Zhao, Jiman. Maximal regularity and Hardy spaces.  Collect. Math.  59  (2008),  no. 1, 103--127. 


\bibitem{B} Bownik, Marcin. Boundedness of operators on Hardy spaces via atomic decompositions.  Proc. Amer. Math. Soc.  133  (2005),  no. 12, 3535--3542.

\bibitem{CD} Coulhon, Thierry; Duong, Xuan Thinh. Maximal regularity and kernel bounds: observations on a theorem by Hieber and Pr\"uss.  Adv. Differential Equations  5  (2000),  no. 1-3, 343--368.

\bibitem{MSV} Meda, Stefano; Sj\"ogren, Peter; Vallarino, Maria. On the $H\sp 1$-$L\sp 1$ boundedness of operators.  Proc. Amer. Math. Soc.  136  (2008),  no. 8, 2921--2931.

\bibitem{M} Meyer, Yves. Ondelettes et op\'erateurs. II. Hermann, Paris, 1990. 

\bibitem {MTW} Meyer, Yves; Taibleson, Mitchell H.; Weiss, Guido.
Some functional analytic properties of the spaces $B\sb q$ generated by blocks.
Indiana Univ. Math. J. 34 (1985), no. 3, 493--515.   




\bibitem{YZ} Yang, Dachun; Zhou, Yuan. A boundedness criterion via atoms for linear operators in Hardy spaces,
to appear in Constr. Approx. DOI 10.1007/s00365-008-9015-1

 \end{thebibliography}
\end{document}